\documentclass{amsart}
\usepackage{graphicx} 

\usepackage[T1]{fontenc}
\usepackage[utf8]{inputenc}
\usepackage[english]{babel}
\usepackage{csquotes}

\usepackage{biblatex}
\renewbibmacro{in:}{%
        \ifentrytype{article}{}{\printtext{\bibstring{in}\intitlepunct}}}

\usepackage[hang,flushmargin]{footmisc} 

\usepackage{amsmath}
\usepackage{amssymb}
\usepackage{amsthm}

\usepackage{abstract}
\addto\captionsenglish{}

\newtheorem{defn}{Definition}[section]
\newtheorem{thm}[defn]{Theorem}
\newtheorem{prop}[defn]{Proposition}
\newtheorem{lemma}[defn]{Lemma}
\newtheorem{cor}[defn]{Corollary}

\newtheorem*{quest}{Question}

\newcommand{\N}{\mathbb{N}}

\newcommand{\R}{\mathbb{R}}

\newcommand{\func}[3]{#1:#2\rightarrow#3}

\DeclareMathOperator{\diam}{diam}
\DeclareMathOperator{\conv}{conv}
\DeclareMathOperator{\sspan}{span}
\DeclareMathOperator{\SO}{SO}

\renewcommand{\epsilon}{\varepsilon}
\renewcommand{\theta}{\vartheta}

\title{On $k$-convex hulls}
\author{Davide Ravasini}
\address{Mathematical Institute, Charles University, Sokolovsk\'{a} 49/83, 186 75 Prague 8, Czech Republic}
\email{davide.ravasini@matfyz.cuni.cz \newline\indent https://orcid.org/0000-0001-6355-5802}

\begin{document}
\maketitle
\let\thefootnote\relax\footnote{\today \newline \indent\emph{2020 Mathematics Subject Classification:} 52A20, 52A23. \newline
}

\begin{abstract}
\noindent \textsc{Abstract}. For every integer $k\geq 2$ and every $R>1$ one can find a dimension $n$ and construct a symmetric convex body $K\subset\R^n$ with $\diam Q_{k-1}(K)\geq R\cdot\diam Q_k(K)$, where $Q_k(K)$ denotes the $k$-convex hull of $K$. The purpose of this short note is to show that this result due to E.\ Kopeck\'{a} is impossible to obtain if one additionally requires that all isometric images of $K$ satisfy the same inequality. To this end, we introduce the dual construction to the $k$-convex hull of $K$, which we call the $k$-cross approximation of $K$. We also prove an infinite-dimensional version of the main result that holds in general Hilbert spaces.
\end{abstract}

\section{Introduction}
A \emph{convex body} is a compact convex subset of $\R^n$ with a nonempty interior. Let $\{e_1,\dots,e_n\}$ be the standard basis of $\R^n$ and, for every $S\subseteq\{1,\dots,n\}$, let us denote by $\func{P_S}{\R^n}{\R^n}$ the orthogonal projection onto the subspace spanned by $\{e_j\,:\,j\in S\}$. For every $k\in\{1,\dots,n\}$, let us call ${[n]}_k$ the family of all subsets $S\subseteq\{1,\dots,n\}$ with cardinality $\# S=k$. For every convex body $K\subset\R^n$, the \emph{$k$-convex hull} of $K$ is defined as
\begin{equation}
    \label{eqn:kconvex}
    Q_k(K)=\bigcap_{S\in{[n]}_k}P_S^{-1}\bigl(P_S(K)\bigr). 
\end{equation}

The reason behind this terminology is easily explained. Given a positive integer $n$ and $k\in\{1,\dots,n\}$, we say that a convex body $K\subset\R^n$ is \emph{$k$-convex} if every $x\in\partial K$ has an outward normal vector with at most $k$ nonzero coordinates. Alternatively, $K$ is $k$-convex if it is the intersection of a family $\mathcal{F}$ of half spaces such that the boundary of every $H\in\mathcal{F}$ is a hyperplane defined by at most $k$ coordinates. It is clear that the intersection of an arbitrary family of $k$-convex bodies is again $k$-convex. This enables us to consider the $k$-convex hull of a convex body $K\subset\R^n$, that is, the intersection of all $k$-convex bodies containing $K$. It is not difficult to check that the $k$-convex hull of $K$ is given precisely by (\ref{eqn:kconvex}).

The $1$-convex hull $Q_1(K)$ is the smallest box aligned with the standard orthonormal system that contains $K$. The family ${\{Q_k(K)\}}_{k=1}^n$ forms a nested sequence of sets which approximates $K$ from above. At the final degree of approximation we have $Q_n(K)=K$. Moreover, if $K$ is symmetric, so are all its $k$-convex hulls. The initial interest in $k$-convex sets lies in the study of contractive sets, that is, the ranges of $1$-Lipschitz retractions. In $\ell_2$ spaces, all closed convex sets are contractive. In $\ell_p$ spaces with $p\in(1,2)\cup(2,\infty)$, the contractive sets are precisely the $2$-convex sets (see \cite{LH, DE, E}). In $\ell_1$, the family of contractive sets contains all $2$-convex sets, but is larger \cite{W}.

If the dimension of the space is allowed to be arbitrarily large, then the diameter of $Q_1(K)$ can be arbitrarily larger than the diameter of $K=Q_n(K)$, as the example of a Euclidean unit ball shows. Surprisingly, this is true not only when one compares $Q_1(K)$ and $K$. In fact, it suffices to consider two successive approximations, as shown by Kopeck\'{a} in \cite[Theorem 2.2]{EK}. We may also measure the diameters involved with respect to any $p$-norm with $p\in[1,\infty)$.

\begin{thm}[Kopeck\'{a}]
    \label{thm:kopecka}
    For every $p\in[1,\infty)$, every integer $k\geq 2$ and every $R>1$ there are an integer $n\geq k$ and a symmetric convex body $K\subset\ell_p^n$ such that
    \[ \frac{\diam Q_{k-1}(K)}{\diam Q_k(K)}\geq R. \]
\end{thm}

However, the construction of $k$-convex hulls is carried out with reference to a single orthonormal basis, which is understood to be the standard basis of $\R^n$. If $K\subset\R^n$ is one of the convex bodies whose existence is predicted by Theorem \ref{thm:kopecka}, it might well be that, for some $U\in\SO(n)$, the body $UK$ fails to satisfy the inequality of Theorem \ref{thm:kopecka}, even though, at least in the Euclidean norm, $K$ and $UK$ are isometric. In a private exchange with the author, E.\ Kopeck\'{a} brought up the following question.
\medskip
\begin{quote}
    Does there exist an absolute constant $C>0$ such that, for every integer $n\geq 2$ and every symmetric convex body $K\subset\R^n$, there is $U\in\SO(n)$ with
    \[ \frac{\diam{Q_1(UK)}}{\diam{Q_2(UK)}}\leq C\,? \]
\end{quote}
\medskip
In this note, we provide an affirmative answer by proving something stronger.
\begin{thm}
    \label{thm:main}
    For every integer $n\geq 2$ and every symmetric convex body $K\subset\R^n$ there exists $U\in\SO(n)$ such that the inclusion
    \[ Q_{k-1}(UK)\subseteq 3Q_k(UK) \]
    holds for every $k\in\{2,\dots,n\}$.
\end{thm} 

As a direct consequence, we see that Theorem \ref{thm:kopecka} fails if $R>3$ and we additionally require that the inequality of the statement be satisfied by all rotations of $K$ simultaneously.
\begin{cor}
    For every positive integer $n$ and for every symmetric convex body $K\subset\R^n$, there is $U\in\SO(n)$ such that the inequality
    \begin{equation}
        \label{eqn:diambound}
        \frac{\diam Q_{k-1}(UK)}{\diam Q_k(UK)}\leq 3
    \end{equation}
    holds for every $k\in\{2,\dots,n\}$, where the diameters may be taken with respect to any norm in $\R^n$.
\end{cor}

Given the independence of the constant $3$ on the dimension of the space in Theorem \ref{thm:main}, it is natural to look for infinite-dimensional generalisations of Theorem \ref{thm:main}. In fact, we will show that a similar result holds in real Hilbert spaces of any density character. Recall that a convex body in a real Hilbert space $H$ is a closed, convex, and bounded set $C\subset H$ with a nonemtpy interior. A complete orthonormal system $\mathbf{e}={\{e_\gamma\}}_{\gamma\in\Gamma}$ over an index set $\Gamma$ is a set of unit vectors in $H$ such that $\sspan(\mathbf{e})$ is dense in $H$ and such that $\langle e_\alpha,e_\beta\rangle=0$ whenever $\alpha\ne\beta$. For every positive integer $k$, we denote by ${[\Gamma]}_k$ the family of all sets $S\subseteq\Gamma$ with $\# S=k$. In line with the preceding notation, given any set $S\subseteq\Gamma$, we denote by $\func{P_S}{H}{H}$ the orthogonal projection onto the subspace $\overline{\sspan}\{e_\gamma\,:\,\gamma\in S\}$. Finally, for every positive integer $k$, let us define the $k$-convex hulls of $C$ with respect to $\mathbf{e}$ as
\[ Q_k(\mathbf{e},C)=\bigcap_{S\in{[\Gamma]}_k}P_S^{-1}\bigl(P_S(C)\bigr). \]

Note that $Q_k(\mathbf{e},C)$ need not be bounded if $H$ is infinite-dimensional. For example, if $B$ is the unit ball of $\ell_2$ and $\mathbf{e}$ is the standard basis of $\ell_2$, then 
\[ Q_1(\mathbf{e},B)=\{x\in\ell_2\,:\,-1\leq x(n)\leq 1 \text{ for every }n\in\N\} \]
and it is easy to see that $Q_1(\mathbf{e},B)$ contains elements of $\ell_2$ of arbitrarily large norm. However, this is not of concern, since the statement of Theorem \ref{thm:main}, which we aim to generalise, also makes sense for unbounded sets.

Let us briefly show that, at least, $Q_k(\mathbf{e},C)$ is a closed set.
\begin{prop}
    Let $C$ be a closed, convex, bounded set with nonempty interior in a Hilbert space $H$ and let $\mathbf{e}$ be an orthonormal system for $H$. Then $Q_k(\mathbf{e},C)$ is closed in the norm topology of $H$ for every positive integer $k$.
\end{prop}
\begin{proof}
Since orthogonal projections are continuous with respect to the weak topology on $H$, it suffices to show that for every weakly continuous map $\func{f}{H}{H}$ the set $f^{-1}(f(C))$ is closed.  Since $C$ is closed, convex, and bounded and $H$ is reflexive, $C$ it is weakly compact. Hence $f(C)$ is weakly compact, thus weakly closed. It follows that $f^{-1}(f(C))$ is also weakly closed and, therefore, closed in the norm topology.
\end{proof}

The theorem we will prove reads as follows.
\begin{thm}
    \label{thm:main2}
    Let $C$ be a symmetric convex body in a Hilbert space $H$. For every $\delta>0$, there exists a complete orthonormal system $\mathbf{e}$ in $H$ such that
    \[ Q_{k-1}(\mathbf{e},C)\subseteq(3+\delta)Q_k(\mathbf{e}, C) \]
    for every integer $k\geq 2$.
\end{thm}

We will present two separate proofs for the finite-dimensional (Theorem \ref{thm:main}) and infinite-dimensional (Theorem \ref{thm:main2}) cases, respectively. In the first step of the proof of Theorem \ref{thm:main}, we will introduce the dual construction to the $k$-convex hull of a symmetric convex body $K$, which we call the $k$-cross approximation of $K$ and denote by $R_k(K)$. In the second step of the proof, we will describe a simple greedy algorithm to find a positively oriented orthonormal basis $\{u_1,\dots,u_n\}$ such that $R_k(UK)\subseteq 3R_{k-1}(UK)$ for every $k\in\{2,\dots,n\}$, where $U\in\SO(n)$ is the isometry of $\R^n$ which maps $u_i$ to $e_i$ for every $i\in\{1,\dots,n\}$. Applying this procedure to the body $K^\circ$, the polar body of $K$, we are finally able to obtain Theorem \ref{thm:main}. The proof of Theorem \ref{thm:main} already contains all the basic ideas needed for the proof of Theorem \ref{thm:main2}. The infinite-dimensional case will differ from the finite-dimensional one only by some technical challenges, namely the necessity to use transfinite induction and the lack of compactness, which is ultimately responsible for the presence of the quantity $\delta$ in the statement of Theorem \ref{thm:main2}.

We briefly recall some notation and terminology that we will use throughout the paper. The reader can find more details in the standard references \cite{gruber, schneider}. The standard Euclidean norm of a vector $x\in\R^n$ will be denoted by $\|x\|$ and the usual scalar product of two vectors $x,y\in\R^n$ by $\langle x,y\rangle$. We will make use of the radial function $\func{\rho_K}{S^{n-1}}{[0,\infty)}$ of a symmetric convex body $K\subseteq\R^n$, given by $\rho_K(u)=\max\{t\geq 0\,:\,tu\in K\}$ for every $u\in S^{n-1}$. We also set
\[ \|K\|=\sup_{x\in K}\|x\|=\sup_{u\in S^{n-1}}\rho_K(u). \]
We denote by $\conv A$ the convex hull of a set $A\subseteq\R^n$, i.e., the smallest convex set that contains $A$. Given two sets $A,B\subseteq\R^n$, their \emph{Minkowski sum} is given by
\[ A+B=\{a+b\,:\,a\in A,\,b\in B\}. \]
Given $x,y\in\R^n$, $[x,y]$ denotes the line segment $\{(1-\lambda)x+\lambda y\,:\,\lambda\in[0,1]\}$. Finally, we denote by $A^\circ$ the \emph{polar} of a set $A\subset\R^n$. This is given by
\[ A^\circ=\{x\in\R^n\,:\,\langle x,y\rangle\leq 1\text{ for every }y\in A\} \]
$A^\circ$ is always closed and convex, and is bounded whenever $0$ is in the interior of $A$. Furthermore, the polar set of a symmetric convex body $K$ is again a symmetric convex body and $K^{\circ\circ}=K$. All these definitions admit straightforward extensions to infinite-dimensional Hilbert spaces. The unit ball and the unit sphere of a general Hilbert space $H$ will be denoted by $B_H$ and $S_H$, respectively. All the Hilbert spaces that we consider will be real.\\

\noindent \textsc{Acknowledgements.} We thank Prof.\ Eva Kopeck\'{a} for being the source of this original, entertaining problem. The author has been supported by the Charles University grant PRIMUS/24/SCI/009.

\section{Cross approximations of a convex body}
\label{sec:finite}
Given an integer $n\geq 2$ and $S\subseteq\{1,\dots,n\}$, let us denote by $F(S)$ the subspace of $\R^n$ spanned by $\{e_j\,:\,j\in S\}$. For every $k\in\{1,\dots,n\}$ and every symmetric convex body $K\subset\R^n$ we set
\begin{equation}
    \label{eqn:cross}
    R_k(K)=\conv\biggl(K\cap\bigcup_{S\in{[n]}_k}F(S)\biggr)
\end{equation}  
and call $R_k(K)$ the \emph{$k$-cross approximation} of $K$. $R_1(K)$ is the largest symmetric cross polytope aligned with the standard orthonormal basis contained in $K$. As in the case of $k$-convex hulls, the family ${\{R_k(K)\}}_{k=1}^n$ forms a nested sequence of convex bodies that approximate $K$ more and more precisely, this time from below, until we get $R_n(K)=K$. 

The duality relation between $k$-convex hulls and $k$-cross approximations is given by the following lemma.
\begin{lemma}
    \label{lemma:cross}
    For every integer $n\geq 2$, every symmetric convex body $K\subset\R^n$ and every $k\in\{1,\dots,n\}$ one has ${Q_k(K)}^\circ=R_k(K^\circ)$.
\end{lemma}
\begin{proof}
    Let us first check that, for every $S\in{[n]}_k$, one has
    \[ {P_S^{-1}\bigl(P_S(K)\bigr)}^\circ=K^\circ\cap F(S). \]
    Clearly, we have $K\subseteq P_S^{-1}(P_S(K))$ and ${F(S)}^\circ={F(S)}^\perp\subseteq P_S^{-1}(P_S(K))$, which implies ${P_S^{-1}(P_S(K))}^\circ\subseteq K^\circ\cap F(S)$. To see the reverse inclusion, let $x\in K^\circ\cap F(S)$ and $y\in P_S^{-1}(P_S(K))$. We may write $y=y_1+y_2$, where $y_1\in K$ and $y_2\in{F(S)}^\perp$. Now we have 
    \[ \langle x,y\rangle =\langle x,y_1\rangle+\langle x,y_2\rangle=\langle x,y_1\rangle\leq 1. \]
    Since $y\in K$ is arbitrary, it follows that $x\in{P_S^{-1}(P_S(K))}^\circ$.
    The statement is now the result of the following simple computations involving polarity.
    \begin{align*} 
        {Q_k(K)}^\circ &={\biggl(\,\bigcap_{S\in{[n]}_k}P_S^{-1}\bigl(P_S(K)\bigr)\biggr)}^\circ \\
        &= \conv\biggl(\,\bigcup_{S\in{[n]}_k}{P_S^{-1}\bigl(P_S(K)\bigr)}^\circ\biggr) \\
        &= \conv\biggl(\,\bigcup_{S\in{[n]}_k}K^\circ\cap F(S)\biggr)=R_k(K^\circ). \qedhere
    \end{align*}
\end{proof}
We also recall that polarity reverses set inclusions and that ${(UK)}^\circ=UK^\circ$ for every symmetric convex body $K\subset\R^n$ and every $U\in\SO(n)$. Thus, we can prove Theorem \ref{thm:main} by showing that for every symmetric convex body $K\subset\R^n$ there is $U\in\SO(n)$ such that $R_k(UK)\subseteq 3R_{k-1}(UK)$ for every $k\in\{1,\dots,n\}$. Let us start with two technical lemmas, which we prove directly in a general Hilbert space.
\begin{lemma}
    \label{lemma:step1}
    Let $H$ be a Hilbert space with $\dim H\geq 2$, let $e\in S_H$, let $\eta\in[0,1)$  and let $L\subset H$ be a symmetric convex body such that $\rho_L(e)\geq(1-\eta)\|L\|$. Then we have
    \[ (1-\eta)L\subseteq L\cap\sspan\{e\}+2L\cap e^\perp. \]
\end{lemma}
\begin{proof}
    If $x\in L$, then $|\langle x,e\rangle|\leq\|L\|$, which means that we can write $x=te+y$, where $|t|\leq\|L\|$ and $y\in e^\perp$. Since $-\rho_L(e)e\in L$, by convexity we have
    \[ -(1-\lambda)\rho_L(e)e+\lambda te+\lambda y\in L \]
    for every $\lambda\in[0,1]$. In particular, if we set 
    \[ \lambda=\frac{\rho_L(e)}{\rho_L(e)+t}, \]
    we obtain 
    \[ \frac{\rho_L(e)}{\rho_L(e)+t}y\in L. \]
    Now we observe that
    \[ \frac{\rho_L(e)}{\rho_L(e)+t}\geq\frac{(1-\eta)\|L\|}{2\|L\|}=\frac{1-\eta}{2}, \]
    which implies that $(1-\eta)y\in 2L$. We also note that
    \[ te=\frac{t}{\rho_L(e)}\rho_L(e)e\in\frac{t}{\rho_L(e)}L\subseteq\frac{\|L\|}{(1-\eta)\|L\|}L=\frac{1}{1-\eta}L, \]
    which implies $(1-\eta)te\in L$. Combining these two observations, we get    \[ (1-\eta)x\in L\cap\sspan\{e\}+2L\cap e^\perp. \] 
    Since $x\in L$ has been arbitrarily chosen, the lemma follows.
\end{proof}
\begin{lemma}
    \label{lemma:step2}
    Let $e$ be a unit vector in a Hilbert space $H$. For every $r>0$ and every symmetric convex body $M\subset e^\perp$ one has
    \[ [-re,re]+2M\subseteq 3\conv\bigl([-re,re]\cup M\bigr). \]
\end{lemma}
\begin{proof}
    Let $\mu$ be the norm on $e^\perp$ generated by $2M$ and let $\nu_1$ and $\nu_2$ be the norms on $H$ generated by the symmetric convex bodies
    \[ [-re,re]+2M\quad\text{and}\quad\conv\bigl([-re,re]\cup M\bigr) \]
    respectively. For every $x\in H$, denote by $p(x)$ its orthogonal projection onto $e^\perp$. Then one has, for every $x\in H$,
    \begin{align*} 
        \nu_2(x) &= \frac{|\langle x,e\rangle|}{r}+2\mu\bigl(p(x)\bigr) \\
        &= 3\biggl(\frac{1}{3}\cdot\frac{|\langle x,e\rangle|}{r}+\frac{2}{3}\cdot\mu\bigl(p(x)\bigr)\biggr) \\
        &\leq3\max\biggl\{\frac{|\langle x,e\rangle|}{r},\,\mu\bigl(p(x)\bigr)\biggr\}=3\nu_1(x).
    \end{align*}
    This inequality for the norms $\nu_1$ and $\nu_2$ proves the claim.
\end{proof}
We can now prove the theorem which, together with Lemma \ref{lemma:cross}, directly implies Theorem \ref{thm:main}.
\begin{thm}
    \label{thm:maincross}
    For every integer $n\geq 2$ and every symmetric convex body $K\subset\R^n$ there is $U\in\SO(n)$ such that $R_k(UK)\subseteq 3R_{k-1}(UK)$ for every $k\in\{2,\dots,n\}$.
\end{thm}
\begin{proof}
    We will construct by induction an orthonormal basis $\{u_1,\dots,u_n\}$ such that the isometry $U\in\SO(n)$ that maps $u_j$ to $e_j$ for every $j$ satisfies the statement. We select $u_1\in S^{n-1}$ such that $\rho_K(u_1)=\max_{u\in S^{n-1}}\rho_K(u)$. Let now $k\in\{2,\dots,n\}$ and assume that $u_j$ has already been defined for $j\in\{1,\dots,k-1\}$. Let $X_k$ be the orthogonal complement of $\sspan\{u_j\,:\,1\leq j\leq k-1\}$ and select $u_k\in S^{n-1}$ such that $\rho_K(u_k)=\max_{u\in S^{n-1}\cap X_k}\rho_K(u)$. If $k=n$, replace $u_n$ with $-u_n$ if necessary to ensure that the bases $\{u_1,\dots,u_n\}$ and $\{e_1,\dots,e_n\}$ share the same orientation. As mentioned earlier, we now pick the unique $U\in\SO(n)$ such that $U(u_j)=e_j$ for every $j\in\{1,\dots,n\}$. Given $k\in\{2,\dots,n\}$, we now want to show that $R_k(UK)\subseteq 3R_{k-1}(UK)$. Using the definition of a $k$-cross approximation given by (\ref{eqn:cross}), it suffices to show that $UK\cap F(S)\subseteq 3R_{k-1}(UK)$ for every $S\in{[n]}_k$. Let then $S\in{[n]}_k$ and set $j_0=\min S$. The goal is to show that $L\subseteq 3R_{k-1}(UK)$, where $L=UK\cap F(S)$. Define further $M=UK\cap\sspan\{e_j\,:\,j\in S,\,j\ne j_0\}$. Observe that, by our choice of $\{u_1,\dots,u_n\}$, we have $\rho_L(e_{j_0})=\max_{u\in S^{n-1}\cap F(S)}\rho_L(u)$, hence we can set $r=\rho_L(e_{j_0})$ and apply Lemma \ref{lemma:step1} (with $\eta=0$) and Lemma \ref{lemma:step2} to $L$ and $M$ respectively to find
    \begin{equation}
        \label{eqn:key}
        L\subseteq[-re_{j_0},re_{j_0}]+2M\subseteq 3\conv\bigl([-re_{j_0},re_{j_0}]\cup M\bigr).
    \end{equation}
    Moreover, we have 
    \[ [-re_{j_0},re_{j_0}]\cup M\subseteq UK\cap\bigcup_{T\in{[n]}_{k-1}}F(T), \]
    which implies $\conv([-re_{j_0},re_{j_0}]\cup M)\subseteq R_{k-1}(UK)$. Together with (\ref{eqn:key}), this concludes the proof.
\end{proof}

\section{The infinite-dimensional case}
In the case of a general Hilbert space, we first have to adjust the definition of the $k$-cross approximation to obtain an infinite-dimensional formulation of Lemma \ref{lemma:cross}. Let $C$ be a symmetric convex body in a Hilbert space $H$ with $\dim H\geq 2$ and let $\mathbf{e}={\{e_\gamma\}}_{\gamma\in\Gamma}$ be a complete orthonormal system for $H$. For every positive integer $k$, the meaning of ${[\Gamma]}_k$ is the same as in the introduction, and, similarly to Section \ref{sec:finite}, we set $F(S)=\overline{\sspan}\{e_\gamma\,:\,\gamma\in S\}$ for every $S\subseteq\Gamma$. For every positive integer $k$, the $k$-cross approximation of $C$ with respect to $\mathbf{e}$ is now defined as
\[ R_k(\mathbf{e},C)=\overline{\conv}\biggl(C\cap\bigcup_{S\in{[\Gamma]}_k}F(S)\biggr). \]
Note that, in general, $R_k(\mathbf{e},C)$ might not be a convex body, in that it might have an empty interior. For example, if $B$ is the unit ball of $\ell_2$ and $\mathbf{e}$ is the standard basis of $\ell_2$, then $R_1(\mathbf{e},B)$ is the unit ball of $\ell_1\subset\ell_2$, which has an empty interior.

The proof of Lemma \ref{lemma:cross} can now be carried out with minor adjustments in the case of a general Hilbert space and yields the following.
\begin{lemma}
    \label{lemma:cross2}
    Let $H$ be a Hilbert space with $\dim H\geq 2$. For every symmetric convex body $C\subset H$, every complete orthonormal system $\mathbf{e}$ and every positive integer $k$ one has ${Q_k(\mathbf{e},C)}^\circ=R_k(\mathbf{e},C^\circ)$.
\end{lemma}
Similarly to how Lemma \ref{lemma:cross} and Theorem \ref{thm:maincross} together imply Theorem \ref{thm:main}, Lemma \ref{lemma:cross2} and the following will directly imply Theorem \ref{thm:main2}.
\begin{thm}
    Let $C$ be a symmetric convex body in a Hilbert space $H$. For every $\delta>0$ there exists an orthonormal system $\mathbf{e}$ in $H$ such that, for every integer $k\geq 2$, one has $R_k(\mathbf{e},C)\subseteq(3+\delta)R_{k-1}(\mathbf{e},C)$.
\end{thm}
\begin{proof}
    Let $c$ be the cardinality of a complete orthonormal system of $H$ and let $\alpha$ be an ordinal number with cardinality greater than $c$. We will use transfinite induction to construct a set $\{e_\beta\}_{\beta<\alpha}$ such that one of its initial segments will be the orthonormal system we are looking for. First, choose $\eta\in(0,1)$ with 
    \[ \frac{3\eta}{(1-\eta)}<\delta. \] 
    Then choose $e_0$ so that $\rho_C(e_0)>(1-\eta)\|C\|$. Now, let $\beta<\alpha$ and assume that $e_\gamma$ has already been determined for every $\gamma<\beta$. Set $X=\overline{\sspan}{\{e_\gamma\}}_{\gamma<\beta}$. If $X=H$, let $e_\beta=0$. If $X\ne H$, let $D=C\cap X^\perp$ and choose $e_\beta\in S_{X^\perp}$ such that $\rho_C(e_\beta)=\rho_D(e_\beta)>(1-\eta)\|D\|$. Note that there must be $\beta<\alpha$ with $e_\beta=0$, otherwise ${\{e_\beta\}}_{\beta<\alpha}$ would be an orthonormal set with cardinality greater than $c$, which is impossible. Define $\beta_0=\min\{\beta<\alpha\,:\,e_\beta=0\}$. We now claim that $\mathbf{e}={\{e_\beta\}}_{\beta<\beta_0}$ is the compete orthonormal system that we wanted to find. Let $k\geq 2$ be an integer, and let $S\subseteq\{\beta\,:\,\beta<\beta_0\}$ be any subset with $\# S=k$. We want to show that $L=C\cap F(S)\subseteq(3+\delta)R_{k-1}(\mathbf{e},C)$. Let $\gamma=\min S$ and $M=C\cap\sspan\{e_\beta\in S\,:\,\beta\in S,\,\beta\ne\gamma\}$. Note that, in the subspace $F(S)$, we have $\rho_L(e_\gamma)>(1-\eta)\|L\|$, which implies, by Lemma \ref{lemma:step1} and \ref{lemma:step2}, that 
    \begin{align*} 
        (1-\eta)L &\subseteq L\cap\sspan\{e_\gamma\}+2M \\
        &= [-\rho_L(e_\gamma)e_\gamma,\rho_L(e_\gamma)e_\gamma]+2M \\
        &\subseteq 3\conv([-\rho_L(e_\gamma)e_\gamma,\rho_L(e_\gamma)e_\gamma]\cup M).
    \end{align*}
    Note that $\conv([-\rho_L(e_\gamma)e_\gamma,\rho_L(e_\gamma)e_\gamma]\cup M)\subseteq R_{k-1}(\mathbf{e},C)$. By our choice of $\eta$, this yields $L\subseteq(3+\delta)R_{k-1}(\mathbf{e},C)$, as we wanted. Since $k$ and $S$ have been picked arbitrarily, the definition of the $k$-cross approximation allows us to conclude that, for every positive integer $k$, we have $R_k(\mathbf{e},C)\subseteq(3+\delta)R_{k-1}(\mathbf{e},C)$, as desired.
\end{proof}

We conclude with a remark regarding an alternative construction of the $k$-convex hull in infinite-dimensional Hilbert spaces. Given a symmetric, convex body $C$ in an infinite-dimensional Hilbert space $H$ and a complete orthonormal system $\mathbf{e}={\{e_\gamma\}}_{\gamma\in\Gamma}\subset H$, instead of projecting $C$ onto finite-dimensional coordinate subspaces, we can project $C$ onto coordinate subspaces of finite codimension. More precisely, we can define, for every $k\in\N$,
\[ Q^k(\mathbf{e},C)=\bigcap_{S\in{[\Gamma]}_k}P_{\Gamma\setminus S}^{-1}\bigl(P_{\Gamma\setminus S}(C)\bigr). \]
However, this is not an interesting notion, as clarified in the next proposition.
\begin{prop}
    Let $C$ be a symmetric convex body in an infinite-dimensional Hilbert space $H$. For every complete orthonormal system $\mathbf{e}\subset H$ and every $k\in\N$ one has $Q^k(\mathbf{e},C)=C$.
\end{prop}
\begin{proof}
    We will work again with cross approximations. For every symmetric convex body $C\subset H$, every $k\in\N$ and every complete orthonormal system $\mathbf{e}={\{e_\gamma\}}_{\gamma\in\Gamma}$, we define
    \[ R^k(\mathbf{e},C)=\overline{\conv}\biggl(C\cap\bigcup_{S\in{[\Gamma]}_k}F(\Gamma\setminus S)\biggr). \]
    It is not difficult to follow the proof of Lemma \ref{lemma:cross} and show that ${Q^k(\mathbf{e},C)}^\circ=R^k(\mathbf{e},C^\circ)$ for every symmetric convex body $C\subset H$, every $k\in\N$ and every complete orthonormal system $\mathbf{e}$, hence the equality of the statement is equivalent to $R^k(\mathbf{e},C)=C$. The inclusion $R^k(\mathbf{e},C)\subseteq C$ is clear and we note that the opposite inclusion is almost trivial if $H$ is nonseparable.  Indeed, in this case $\Gamma$ is uncountable, but, for every $x\in H$, the set $T_x=\{\gamma\in\Gamma\,:\,\langle x,e_\gamma\rangle\ne 0\}$ is at most countable. Therefore, for every $k\in\N$ there is a subset $S\subset\Gamma$ with $\# S=k$ and $T_x\subset\Gamma\setminus S$, that is, $x\in F(\Gamma\setminus S)$. The case of a countable orthonormal system $\mathbf{e}={\{e_n\}}_{n\in\N}$ remains to be discussed. Since $C$ is a symmetric convex body, there is $\rho>0$ such that $\rho B_H\subseteq C$. Let $x\in C$. We want to find a sequence ${(x_j)}_{j=1}^\infty$ such that $x_j\in R^k(\mathbf{e},C)$ for every $j$ and $x_j\to x$ as $j\to\infty$. It will then follow that $x\in R^k(\mathbf{e},C)$, as $R^k(\mathbf{e},C)$ is closed. To this end, for every positive integer $j$ find $n$ such that 
    \[ \sum_{i=n}^\infty{\langle x,e_i\rangle}^2<\frac{1}{j^2} \]
    Then, define $y\in H$ as $y=\sum_{i\in\N}\alpha(i)e_i$, where
    \[ \alpha(i)=\left\{\begin{array}{ll}
                        0 & \text{ if }0\leq i<n \\
                        -j\rho \langle x,e_i\rangle & \text{ if }i\geq n 
                        \end{array}\right. \]
    Note that $\|y\|<\rho$, which implies $y\in C$. We now set
    \[ x_j=\frac{j\rho}{j\rho+1}x+\frac{1}{j\rho+1}y \]
    and observe that $x_j\in C$ by convexity. Moreover, $\langle x_j,e_i\rangle=0$ for every $i\geq n$, which means that $x_j\in F(\N\setminus S)$ for every $S\subset\N$ with $S\cap\{0,\dots,n-1\}=\varnothing$. Hence, we surely have $x_j\in R^k(\mathbf{e},C)$. Finally, one can compute that
    \begin{align*}
        \|x-x_j\|=\frac{1}{j\rho+1}\|x-y\|\leq\frac{\|x\|+\|y\|}{j\rho+1}<\frac{\|x\|+\rho}{j\rho+1}\to 0
    \end{align*}
    as $j\to\infty$. Therefore, the sequence ${(x_j)}_{j=1}^\infty$ converges to $x$, as desired.
\end{proof}
\printbibliography

\end{document}

In the first step of the proof, we introduce the dual construction to the $k$-convex hull of a symmetric convex body $K$, which we call the $k$-cross approximation of $K$ and denote by $R_k(K)$. In the second step of the proof, we describe a simple greedy algorithm to find a positively oriented orthonormal basis $\{u_1,\dots,u_n\}$ such that $R_k(UK)\subseteq 3R_{k-1}(UK)$ for every $k\in\{2,\dots,n\}$, where $U\in\SO(n)$ is the isometry of $\R^n$ which maps $u_i$ to $e_i$ for every $i\in\{1,\dots,n\}$. Applying this procedure to the body $K^\circ$, the polar body of $K$, we are finally able to obtain Theorem \ref{thm:main}.

We briefly recall some notation and terminology that we will use throughout the paper. The reader can find more details in the standard references \cite{gruber, schneider}. We will make use of the radial function $\func{\rho_K}{S^{n-1}}{[0,\infty)}$ of a symmetric convex body $K\subseteq\R^n$, given by $\rho_K(u)=\max\{t\geq 0\,:\,tu\in K\}$ for every $u\in S^{n-1}$. We denote by $\conv A$ the convex hull of a set $A\subseteq\R^n$, i.e.\, the smallest convex set that contains $A$. Given two sets $A,B\subseteq\R^n$, their \emph{Minkowski sum} is given by
\[ A+B=\{a+b\,:\,a\in A,\,b\in B\}. \]
Given $x,y\in\R^n$, $[x,y]$ denotes the line segment $\{(1-\lambda)x+\lambda y\,:\,\lambda\in[0,1]\}$. Finally, we denote by $\langle x,y\rangle$ the usual scalar product of two vectors $x,y\in\R^n$ and by $A^\circ$ the \emph{polar} of a set $A\subset\R^n$. This is given by
\[ A^\circ=\{x\in\R^n\,:\,\langle x,y\rangle\leq 1\text{ for every }y\in A\} \]
$A^\circ$ is always closed and convex, and it is compact whenever $0$ lies in the interior of $A$. Furthermore, the polar set of a symmetric convex body $K$ is again a symmetric convex body and $K^{\circ\circ}=K$. \\

\noindent \textsc{Acknowledgements.} We thank Prof.\ Eva Kopeck\'{a} for being the source of this original, entertaining problem. The author has been supported by the Charles University grant PRIMUS/24/SCI/009.

\section{Cross approximations of a convex body}
Given an integer $n\geq 2$ and $S\subseteq\{1,\dots,n\}$, let us denote by $F(S)$ the subspace of $\R^n$ spanned by $\{e_j\,:\,j\in S\}$. For every $k\in\{1,\dots,n\}$ and every symmetric convex body $K\subset\R^n$ we set
\begin{equation}
    \label{eqn:cross}
    R_k(K)=\conv\biggl(K\cap\bigcup_{S\in{[n]}_k}F(S)\biggr)
\end{equation}  
and call $R_k(K)$ the \emph{$k$-cross approximation} of $K$. $R_1(K)$ is the largest symmetric cross polytope aligned with the standard orthonormal basis contained in $K$. As in the case of $k$-convex hulls, the family ${\{R_k(K)\}}_{k=1}^n$ forms a nested sequence of convex bodies that approximate $K$ more and more precisely, this time from below, until we get $R_n(K)=K$. 

The duality relation between $k$-convex hulls and $k$-cross approximations is given by the following lemma.
\begin{lemma}
    \label{lemma:cross}
    For every integer $n\geq 2$, every symmetric convex body $K\subset\R^n$ and every $k\in\{1,\dots,n\}$ one has ${Q_k(K)}^\circ=R_k(K^\circ)$.
\end{lemma}
\begin{proof}
    Let us first check that, for every $S\in{[n]}_k$, one has
    \[ {P_S^{-1}\bigl(P_S(K)\bigr)}^\circ=K^\circ\cap F(S). \]
    Clearly, we have $K\subseteq P_S^{-1}(P_S(K))$ and ${F(S)}^\circ={F(S)}^\perp\subseteq P_S^{-1}(P_S(K))$, which implies ${P_S^{-1}(P_S(K))}^\circ\subseteq K^\circ\cap F(S)$. To see the reverse inclusion, let $x\in K^\circ\cap F(S)$ and $y\in P_S^{-1}(P_S(K))$. We may write $y=y_1+y_2$, where $y_1\in K$ and $y_2\in{F(S)}^\perp$. Now we have 
    \[ \langle x,y\rangle =\langle x,y_1\rangle+\langle x,y_2\rangle=\langle x,y_1\rangle\leq 1. \]
    Since $y\in K$ is arbitrary, it follows that $x\in{P_S^{-1}(P_S(K))}^\circ$.
    The statement is now the result of the following simple computations involving polarity.
    \begin{align*} 
        {Q_k(K)}^\circ &={\biggl(\,\bigcap_{S\in{[n]}_k}P_S^{-1}\bigl(P_S(K)\bigr)\biggr)}^\circ \\
        &= \conv\biggl(\,\bigcup_{S\in{[n]}_k}{P_S^{-1}\bigl(P_S(K)\bigr)}^\circ\biggr) \\
        &= \conv\biggl(\,\bigcup_{S\in{[n]}_k}K^\circ\cap F(S)\biggr)=R_k(K^\circ). \qedhere
    \end{align*}
\end{proof}
We also recall that polarity reverses set inclusions and that ${(UK)}^\circ=UK^\circ$ for every symmetric convex body $K\subset\R^n$ and every $U\in\SO(n)$. Thus, we can prove Theorem \ref{thm:main} by showing that for every symmetric convex body $K\subset\R^n$ there is $U\in\SO(n)$ such that $R_k(UK)\subseteq 3R_{k-1}(UK)$ for every $k\in\{1,\dots,n\}$. Let us start with two technical lemmas.
\begin{lemma}
    \label{lemma:step1}
    Let $k\geq 2$ be an integer and let $L\subset\R^k$ be a symmetric convex body such that $\rho_L(e_k)=\max_{u\in S^{n-1}}\rho_L(u)$, where $\{e_1,\dots,e_k\}$ is the standard basis of $\R^k$. Then we have
    \[ L\subseteq L\cap\sspan\{e_k\}+2L\cap\sspan\{e_j\,:\,1\leq j\leq k-1\}. \]
\end{lemma}
\begin{proof}
    Let $x\in L$. Since $\rho_L$ assumes its largest value at $e_k$ and $-e_k$, we can write $x=te_k+y$ with $|t|\leq\rho_L(e_k)$ and $y\in\sspan\{e_j\,:\,1\leq j\leq k-1\}$. We want to show that $y\in 2L$. We may assume by the symmetry of $L$ that $t\geq 0$. By convexity, the line segment $[-\rho_L(e_k)e_k,x]$ lies entirely in $L$. More explicitly, for every $\lambda\in[0,1]$ we have 
    \[ -(1-\lambda)\rho_L(e_k)e_k+\lambda te_k+\lambda y\in L. \]
    In particular, if we set
    \[ \lambda=\frac{\rho_L(e_k)}{\rho_L(e_K)+t}, \]
    we have $\lambda>1/2$ and $\lambda y\in L$. This implies $y\in 2L$, as we wanted.
\end{proof}
\begin{lemma}
    \label{lemma:step2}
    Let $\{e_1,\dots,e_k\}$ be the standard basis of $\R^k$. For every $r>0$ and every symmetric compact convex set $M\subset\sspan\{e_j\,:\,1\leq j\leq k-1\}$ one has
    \[ [-re_k,re_k]+2M\subseteq 3\conv\bigl([-re_k,re_k]\cup M\bigr). \]
\end{lemma}
\begin{proof}
    Let $\mu$ be the norm on $\sspan\{e_j\,:\,1\leq j\leq k-1\}$ generated by $2M$ and let $\nu_1$ and $\nu_2$ be the norms on $\R^k$ generated by the symmetric convex bodies
    \[ [-re_k,re_k]+2M\quad\text{and}\quad\conv\bigl([-re_k,re_k]\cup M\bigr) \]
    respectively. For every $x\in\R^k$, denote by $p(x)$ its orthogonal projection onto $\sspan\{e_j\,:\,1\leq j\leq k-1\}$. Then one has, for every $x\in\R^k$,
    \begin{align*} 
        \nu_2(x) &= \frac{|\langle x,e_k\rangle|}{r}+2\mu\bigl(p(x)\bigr) \\
        &= 3\biggl(\frac{1}{3}\cdot\frac{|\langle x,e_k\rangle|}{r}+\frac{2}{3}\cdot\mu\bigl(p(x)\bigr)\biggr) \\
        &\leq3\max\biggl\{\frac{|\langle x,e_k\rangle|}{r},\,\mu\bigl(p(x)\bigr)\biggr\}=3\nu_1(x).
    \end{align*}
    This inequality for the norms $\nu_1$ and $\nu_2$ proves the claim.
\end{proof}
We can now prove the theorem which, together with Lemma \ref{lemma:cross}, directly implies Theorem \ref{thm:main}.
\begin{thm}
    For every integer $n\geq 2$ and every symmetric convex body $K\subset\R^n$ there is $U\in\SO(n)$ such that $R_k(UK)\subseteq 3R_{k-1}(UK)$ for every $k\in\{2,\dots,n\}$.
\end{thm}
\begin{proof}
    We will construct by induction an orthonormal basis $\{u_1,\dots,u_n\}$ such that the isometry $U\in\SO(n)$ that maps $u_j$ to $e_j$ for every $j$ satisfies the statement. We select $u_n\in S^{n-1}$ such that $\rho_K(u_n)=\max_{u\in S^{n-1}}\rho_K(u)$. Let now $k\in\{1,\dots,n-1\}$ and assume that $u_{n-j}$ has already been defined for $j\in\{0,\dots,k-1\}$. Let $X_k$ be the orthogonal complement of $\sspan\{u_{n-j}\,:\,0\leq j\leq k-1\}$ and select $u_{n-k}\in S^{n-1}$ such that $\rho_K(u_{n-k})=\max_{u\in S^{n-1}\cap X_k}\rho_K(u)$. If $k=n-1$, replace $u_1$ with $-u_1$ if necessary to ensure that the bases $\{u_1,\dots,u_n\}$ and $\{e_1,\dots,e_n\}$ share the same orientation. As mentioned earlier, we now pick the unique $U\in\SO(n)$ such that $U(u_j)=e_j$ for every $j\in\{1,\dots,n\}$. Given $k\in\{2,\dots,n\}$, we now want to show that $R_k(UK)\subseteq 3R_{k-1}(UK)$. Using the definition of a $k$-cross approximation given by (\ref{eqn:cross}), it suffices to show that $UK\cap F(S)\subseteq 3R_{k-1}(UK)$ for every $S\in{[n]}_k$. Let then $S\in{[n]}_k$ and set $j_0=\max S$. The goal is to show that $L\subseteq 3R_{k-1}(UK)$, where $L=UK\cap F(S)$. Define further $M=UK\cap\sspan\{e_j\,:\,j\in S,\,j\ne j_0\}$. Observe that, by our choice of $\{u_1,\dots,u_n\}$, we have $\rho_L(e_{j_0})=\max_{u\in S^{n-1}\cap F(S)}\rho_L(u)$, hence we can set $r=\rho_L(e_{j_0})$ and apply Lemma \ref{lemma:step1} and Lemma \ref{lemma:step2} to $L$ and $M$ respectively to find
    \begin{equation}
        \label{eqn:key}
        L\subseteq[-re_{j_0},re_{j_0}]+2M\subseteq 3\conv\bigl([-re_{j_0},re_{j_0}]\cup M\bigr).
    \end{equation}
    Moreover, we have 
    \[ [-re_{j_0},re_{j_0}]\cup M\subseteq UK\cap\bigcup_{T\in{[n]}_{k-1}}F(T), \]
    which implies $\conv([-re_{j_0},re_{j_0}]\cup M)\subseteq R_{k-1}(UK)$. Together with (\ref{eqn:key}), this concludes the proof.
\end{proof}

\section{The infinite-dimensional case}

\section{On the optimality of the result}
To investigate the optimality of the bound provided by (\ref{eqn:diambound}), let us consider the Euclidean unit ball $B_n\subseteq\R^n$. In this case, no isometry $U\in\SO(n)$ affects the ratio on the left side of (\ref{eqn:diambound}).
\begin{prop}
    For every positive integer $n$ and every $k\in\{2,\dots,n\}$ one has
    \[ Q_{k-1}(B_n)\subseteq\sqrt{\frac{k}{k-1}}Q_k(B_n). \]
\end{prop}
\begin{proof}
    For every $k\in\{1,\dots,n\}$, the explicit form of $Q_k(B_n)$ is given by
    \begin{equation}
        \label{eqn:ballhulls}
        Q_k(B_n)=\biggl\{(x_1,\dots,x_n)\in\R^n\,:\,\sum_{j\in S}x_j^2\leq 1 \text{ for every }S\in{[n]}_k\biggr\}
    \end{equation} 
     Let $k\in\{2,\dots,n\}$. For every $x\in Q_{k-1}(B)$ and every $S\in{[n]}_k$ one has
    \[ \sum_{j\in S}x_j^2=\frac{1}{k-1}\sum_{j\in S}(k-1)x_j^2=\frac{1}{k-1}\sum_{\substack{T\subset S \\ \# T=k-1}}\sum_{j\in T}x_j^2\leq\frac{k}{k-1}. \]
    This proves the desired inclusion.
\end{proof}
The above proposition allows us to infer that for every integer $n\geq 2$ and every $k\in\{2,\dots,n\}$ we have
\[ \frac{\diam Q_{k-1}(B_n)}{\diam Q_k(B_n)}\leq\sqrt{\frac{k}{k-1}}, \]
where the diameters can again be taken with respect to any norm. This ratio is much better than the $3$ that appears in ($\ref{eqn:diambound}$), since it is a function of $k$ which approaches $1$ as $k\to\infty$. In the case of the Euclidean norm, one cannot do better.
\begin{prop}
    For every positive integer $n$ and every $k\in\{1,\dots,n\}$ one has
    \[ \diam Q_k(B_n)=2\sqrt{\frac{n}{k}}, \]
    where the diameter is taken with respect to the Euclidean norm. In particular,
    \[ \frac{\diam Q_{k-1}(UB_n)}{\diam Q_k(UB_n)}=\sqrt{\frac{k}{k-1}} \]
    for every $U\in\SO(n)$.
\end{prop}
\begin{proof}
    Using again the explicit form of $Q_k(B_n)$ given in (\ref{eqn:ballhulls}) one can easily check that the point
    \[ x=\frac{1}{\sqrt{k}}\sum_{j=1}^ne_j \]
    lies in $Q_k(B_n)$, therefore
    \[ \diam Q_k(B_n)\geq 2|x|=2\sqrt{\frac{n}{k}}. \]
    We now have to show that, for every $y\in Q_k(B_n)$, one has $|y|\leq|x|$. If we write $y_j^2=k^{-1}+\alpha_j$ for every $j$, then we must have $\sum_{j\in S}\alpha_j\leq 0$ for every $S\in{[n]}_k$. Hence
    \begin{align*}
        \binom{n-1}{k-1}{|y|}^2 &= \binom{n-1}{k-1}\sum_{j=1}^n(k^{-1}+\alpha_j) \\
        &= \binom{n-1}{k-1}\frac{n}{k}+\sum_{S\in{[n]}_k}\sum_{j\in S}\alpha_j \\
        &\leq\binom{n-1}{k-1}\frac{n}{k},
    \end{align*} 
    which gives us the desired inequality.
\end{proof}
These considerations lead to the following open question.
\begin{quest}
    Is there a function $\func{f}{\N}{[1,\infty)}$ such that $f(k)\to 1$ as $k\to\infty$ and such that, for every integer $n\geq 2$, every $k\in\{2,\dots,n\}$ and every symmetric convex body $K\subset\R^n$ there is $U\in\SO(n)$ with 
    \[ \frac{\diam Q_{k-1}(UK)}{\diam Q_k(UK)}\leq f(k)\,? \]
\end{quest}
\noindent One might even hope for a stronger result.
\begin{quest}
    Is there a function $\func{f}{\N}{[1,\infty)}$ such that $f(k)\to 1$ as $k\to\infty$ and such that, for every integer $n\geq 2$ and every symmetric convex body $K\subset\R^n$ there is $U\in\SO(n)$ with 
    \[ \frac{\diam Q_{k-1}(UK)}{\diam Q_k(UK)}\leq f(k)\,? \]
\end{quest}